\documentclass{article}

\newtheorem{theorem}{Theorem}[section]
\newtheorem{lemma}[theorem]{Lemma}

\newcommand\prob{{\mbox Prob}}

\newcommand\fphp{{(\neg \mbox{PHP}_n)}}

\newcommand\fp{{{\bf F}_p}}
\newcommand\on{{\Omega_{e, n, r}}}

\newcommand\sne{{S_{n,e}}}

\newcommand\sfn{S^{\cal F}_{n,e}}

\newcommand\eer{{Error^{\cal F}_{n,e}}}
\newcommand\eerr{{Error^{\cal F}_{n,r}}}
\newcommand\sfnr{S^{\cal F}_{n,r}}
\newcommand\om{{\Omega_{{\cal F}_n,r}}}

\newcommand\qed{\begin{flushright} {\bf q.e.d.} \end{flushright} }
\newcommand\prf{\noindent {\bf Proof :}  }

\begin{document}

\title{A reduction of proof complexity to computational complexity
for $AC^0[p]$ Frege systems}

\author{Jan Kraj\'{\i}\v{c}ek}

\date{Faculty of Mathematics and Physics\\
Charles University in Prague}

\maketitle

\begin{abstract}
We give a general reduction of lengths-of-proofs lower bounds for
constant depth Frege systems in DeMorgan language augmented by
a connective counting modulo a prime $p$
(the so called $AC^0[p]$ Frege systems)
to computational complexity
lower bounds for search tasks involving search trees branching upon
values of maps on the vector space of
low degree polynomials over $\fp$.
\end{abstract}

In 1988 Ajtai \cite{Ajt88} proved that the unsatisfiable set 
$\fphp$
of propositional formulas 
$$
\bigvee_{j \in [n]} p_{ij}\ 
\mbox{ and } \ 
\ \neg p_{i_1 j} \vee \neg p_{i_2 j}\ 
\mbox{ and } \ 
\ \neg p_{i j_1} \vee \neg p_{i j_2}
$$
for all $i \in [n+1] = \{1, \dots, n+1\}$ , all
$i_1 \neq i_2 \in [n+1], j \in [n]$, and all
$i \in [n+1], j_1 \neq j_2 \in [n]$ respectively,
expressing the failure of the pigeonhole principle (PHP), has for no $d \geq 1$
a polynomial size refutation in a Frege proof system operating
only with DeMorgan formulas of depth at most $d$. 
Subsequently  Kraj\'{\i}\v{c}ek \cite{Kra-lower} established an
exponential lower bound for these so called $AC^0$ Frege proof systems (for different formulas)
and  Kraj\'{\i}\v{c}ek, Pudl\' ak and Woods \cite{KPW} and
Pitassi, Beame and Impagliazzo \cite{PBI} improved independently (and
announced jointly in \cite{BIKPPW}) Ajtai's bound
for PHP to exponential.

All these papers employ some adaptation of the random restriction method
that has been so successfully applied earlier in circuit complexity
(cf. \cite{Ajt83,FSS, Yao,Has}). Razborov \cite{Raz87}
invented already in 1987
an elegant method, simplified and generalized by Smolensky \cite{Smo},
to prove lower bounds even for $AC^0[p]$ circuits, $p$ a prime.
Thus immediately after the lower bounds for $AC^0$ Frege systems were
proved researchers attempted to adapt the Razborov-Smolensky method
to proof complexity and to prove lower bounds also for $AC^0[p]$ Frege systems.

This turned out to be rather elusive and no lower bounds for the systems 
were proved, although some related results were obtained.
Ajtai \cite{Ajt90,Ajt94a,Ajt94b}, 
Beame et.al.\cite{BIKPP} and Buss et.al.\cite{BIKPRS}
proved lower bounds for $AC^0$ Frege systems in DeMorgan language
augmented by the so called modular counting principles as extra axioms
(via degree lower bounds for the Nullstellensatz proof system 
in \cite{BIKPP,BIKPRS}),
Razborov \cite{Raz98} proved $n/2$ degree lower bound for  
refutations of $\fphp$ in polynomial calculus PC of
Clegg, Edmonds and Impagliazzo \cite{CEI}, and Kraj\'{\i}\v cek \cite{Kra-degree}
used methods of Ajtai \cite{Ajt94a,Ajt94b} to prove $\Omega(\log \log n)$
degree lower bound for PC proofs of the counting principles.
Kraj\'{\i}\v cek \cite{Kra-mfcs} \index{Kra} proved
an exponential lower bound for a subsystem of an
$AC^0[p]$ Frege system that extends both constant depth
Frege systems and polynomial calculus, 
Maciel and Pitassi \cite{MP1} demonstrated a quasi-polynomial simulation
of $AC^0[p]$ proof systems by a proof system operating with depth $3$
threshold formula, 
Impagliazzo and Segerlind \cite{IS}
proved that $AC^0$ Frege systems with counting
axioms modulo a prime $p$ do not polynomially
simulate polynomial calculus over ${\bf F}_p$, and
recently Buss, Kolodziejczyk and Zdanowski \cite{BKZ}
proved that an $AC^0[p]$ Frege system of any fixed depth can be
quasi-polynomially simulated by the depth $3$ $AC^0[p]$ system.
Also, Buss et.al.\cite{BIKPRS} showed that the $AC^0[p]$ Frege systems
are polynomially equivalent to the Nullstellensatz proof system
of Beame et.al.\cite{BIKPP} augmented by the so called extension axioms
formalizing in a sense the Razborov-Smolensky method.

\bigskip

In this paper we reduce the task to prove 
a lengths-of-proofs lower bound 
for $AC^0[p]$ Frege systems 
to the task to establish a 
computational hardness of a specific computational task.
The task is a search task and it is solved by
trees branching upon
values of maps on the vector space of
low degree polynomials over $\fp$. The hardness statements to which
lower bounds are reduced say that every tree
of small depth and using small degree polynomials 
succeeds only with an exponentially small probability.

Maciel and Pitassi \cite{MP2} formulated such a reduction of proof
complexity to computational complexity (and the implied 
conditional lower bound). However, in their construction they needed
to redefine the proof systems (allowing arbitrary formulas 
with $MOD_{p,i}$ connectives and
restricting only cut-formulas to constant depth) and the
hard examples whose short proofs yield a computational information are 
not $AC^0[p]$ formulas. In particular, their reduction 
does not seem to 
yield anything for the originally defined $AC^0[p]$ Frege systems
(see Section \ref{1}).

The paper is organized as follows. 
In Section \ref{1} we recall the definition of the proof systems. 
In Sections \ref{2} - \ref{5} we reduce
the lower bounds to the task to show the existence of 
winning strategies for a certain game. 
This is reduced further in Section \ref{6} to the task
to show that search trees of small depth that branch upon values
of maps on the vector space 
of low degree polynomials over $\fp$ cannot solve a certain specific
computational task. 

More background on proof complexity can be found in \cite{kniha}
or \cite{Pud-survey}, the problem (and a relevant background)
to prove the lower bound
for the systems is discussed in detail also in \cite[Chpt.22]{k2}.

\section{$AC^0[p]$ Frege proof systems}
\label{1}

We will work with a sequent calculus style proof system in a language
with connectives $\neg$, unbounded arity $\bigvee$ and unbounded arity
connectives $MOD_{p,i}$ for $p$ a prime and $i = 0, \dots, p-1$.
The intended meaning of the formula
$MOD_{p,i}(y_1, \dots, y_k)$ is that
$\sum_{i} y_i \equiv i\ (\mbox{mod } p)$.
The proof system has the usual structural rules (weakening, contraction
and exchange), the cut rule, the left and the right $\neg$ introduction rules
and two introduction rules for $\bigvee$ 
modified for the unbounded arity; the
$\bigvee : left$ rule
$$
\frac{\varphi_1, \Gamma \rightarrow \Delta\ \ \ \ 
\varphi_2, \Gamma \rightarrow \Delta\ \ \ \dots\ \  \ 
\varphi_t, \Gamma \rightarrow \Delta}{
\bigvee_{i \le t}\varphi_i, \Gamma \rightarrow \Delta}
$$ 
and the
$\bigvee : right$ rule
$$
\frac{\Gamma \rightarrow \Delta, \varphi_j}{\Gamma \rightarrow \Delta, \bigvee_{i \le t}\varphi_i}
$$ 
any $j \le t$.
There are no rules concerning the $MOD_{p,i}$ connectives but there are
new $\mbox{{\bf MOD}}_p$-{\bf axioms} (we follow \cite[Sec.12.6]{kniha}):
\begin{itemize}
\item $MOD_{p,0}(\emptyset)$

\item $
\neg MOD_{p,i}(\emptyset)\ ,\ \mbox{for}\ i=1,\dots, p-1
$

\item $
MOD_{p,i}(\Gamma,\phi) \ \equiv\ [(MOD_{p,i}(\Gamma) \wedge
\neg \phi) \vee (MOD_{p,i-1}(\Gamma) \wedge \phi)]
$

for $i = 0, \dots, p-1$, where $i-1$ means $i-1$ modulo $p$,
and where $\Gamma$ stands for a sequence (possibly empty)
of formulas.
\end{itemize}
The depth of the formula is the maximal number of alternations of the
connectives; in particular, formulas from $\fphp$ have depth
$1$ and $2$ respectively.
We have not included among the connectives the conjunction $\bigwedge$;
this is in order to decrease the number of cases
one needs to consider in the constructions later on. Note that the need
to express $\bigwedge$ using $\neg$ and $\bigvee$ may increase the depth
of $AC^0$ formulas comparing to how it is usually counted. But as
we are aiming at lower bounds for all depths this is irrelevant.

We shall denote the proof system $LK(MOD_p)$ and its depth $d$ subsystem
(operating only with formulas of depth at most $d$) $LK_d(MOD_p)$.
It is well-known that this system is polynomially equivalent to constant depth
Frege systems with $MOD_{p,i}$ connectives
(or to Tait style system as in \cite{BKZ}) and in the mutual simulation
the depth increases only by a constant as the systems have the same language
(cf.\cite{kniha}). 
The size of a formula or of a proof is the total number of symbols in it.

\section{From a proof to a game with formulas}
\label{2}

In this section and in the next one we define certain games
using the specific case of the PHP as an example. This is in order
no to burden the presentation right at the beginning
with a technical discussion of the 
form of formulas we allow. As it is shown in Section \ref{6}
this is without a loss of generality and, in fact, motivates
the general formulation there.

Consider the following game $G(d,n,t)$ played between two players,
Prover and Liar. At every round Prover asks a question which Liar must
answer. Allowed questions are:
\begin{enumerate}

\item[(P1)] What is the truth-value of $\varphi$?

\item[(P2)] If Liar already gave a truth-value to $\varphi = \bigvee_{i \le u} \varphi_i$,
Prover can ask as follows:

\begin{enumerate}
\item If Liar answered {\sf false} then Prover can ask an extra question
about the truth value of any one of $\varphi_j$, $j \le u$.

\item If Liar answered {\sf true} then Prover can request that Liar witnesses
his answer by giving a $j \le u$ and stating that $\varphi_j$ is
{\sf true}.

\end{enumerate}
\end{enumerate}
All formulas asked by Prover 
are built from the variables of $\fphp$,
and must have the depth at most $d$
and the size at most $2^t$.
The Liar's answers must obey the following rules:
\begin{enumerate}

\item[(L0)] When asked about a formula he already gave a truth value to in an earlier
round Liar must give the same answer. 

\item[(L1)] He must give $\varphi$ and $\neg \varphi$
opposite truth values.

\item[(L2)] If asked according to (P2a) about $\varphi_j$ he must give value 
{\sf false}.
If asked according to (P2b) he must give value
{\sf true} also to
some $\varphi_j$ with  $j \le u$.  

\item[(L3)] If asked about any $MOD_p$-axiom he must say {\sf true}.

\item[(L4)] If asked about any formula from $\fphp$ he must say {\sf true}.
\end{enumerate}

The game runs for $t$ rounds of questions and Liar wins if he can always
answer while obeying the rules. Otherwise Prover wins.

\begin{lemma} \label{2.1}
For any $d \geq 2$, $n \geq 1$ and $s \geq 1$. If there is a size $s$ $LK_d(MOD_p)$
refutation of $\fphp$ then Prover has a winning strategy 
for game 
$$G(d + O(1), n, O(\log s))\ .$$
\end{lemma}

\prf

It is well-known that LK-proofs (or Frege proofs) can be put into a form of balanced
tree with only a polynomial increase in size and a constant increase in the depth
(cf. \cite{Kra-lower,kniha}). In particular, the hypothesis of the lemma
implies that there is a size $s^{O(1)}$ refutation $\pi$ of $\fphp$
in $LK_{d+O(1)}(MOD_p)$ that is in a form of tree whose depth is 
$O(\log s)$. 

The Prover will attempt - by asking Liar suitable questions -
to built a path of sequents $Z_1, Z_2, \dots$ in $\pi$ such that
\begin{itemize}

\item $Z_1$ is the end-sequent of $\pi$, i.e. the empty sequent.

\item $Z_{i+1}$ is one of the hypothesis of the inference yielding $Z_i$.

\item If $Z_i$ is $\Gamma \rightarrow \Delta$ then Prover asked all formulas in
$\Gamma, \Delta$ and Liar asserted that all formulas in $\Gamma$ are true
and all formulas in $\Delta$ are false.

\end{itemize}
Assume $Z_1, \dots, Z_i$ has been constructed. Next Prover's move
depends on the type of inference yielding $Z_i$:
\begin{itemize}

\item Structural rules: Prover asks no questions and just takes for $Z_{i+1}$
the hypothesis of the inference.

\item Cut rule: Prover asks about the truth value of the cut formula, say $\varphi$,
and if Liar asserts it to be true, Prover takes for $Z_{i+1}$ the hypothesis of 
the inference having $\varphi$ in the antecedent, otherwise it takes the hypothesis
with $\varphi$ in the succedent.

\item A $\neg$ introduction rule: if $\neg \varphi$ was the formula introduced,
Prover asks $\varphi$ and takes for $Z_{i+1}$ the unique hypothesis of the inference.

\item The $\bigvee : right$ introduction rule: if the principal formula was 
$\varphi = \bigvee_{i \le u}\varphi_i$ and the minor formula $\varphi_j$
Prover already asked $\varphi$ in an earlier round and got answer {\sf false}.
He now asks $\varphi_j$ and 
takes for $Z_{i+1}$ the unique hypothesis of the inference.

\item The $\bigvee : left$ introduction rule: 
if the principal formula was 
$\varphi = \bigvee_{i \le u}\varphi_i$ Prover already asked $\varphi$ in 
an earlier round and got answer {\sf true}. She now asks Liar to witness
this answer by some $\varphi_j$ and then 
takes for $Z_{i+1}$ the hypothesis with the minor formula $\varphi_j$ in the antecedent.

\end{itemize}
This process either causes Liar to lose or otherwise arrives at an initial 
sequent which Liar's answers claim to be false. But that contradicts
rules (L1), (L3) or (L4). 

\qed

Shallow tree-like refutations of a set of axioms can be used as search trees
finding an axiom false under a given assignment: the Liar answers the truth 
values determined by the assignment (see e.g. the use of such trees
in \cite{Kra-lower,kniha}). 
It was an important insight
of Buss and Pudl\' ak \cite{BusPud} that when Liars are allowed not to follow an
assignment but are only required to be logically consistent then the 
minimal length of Prover's winning strategy
characterizes the minimal depth of a tree-like refutation
(a form of a statement opposite to the lemma also holds as pointed out
in \cite{BusPud} in the context of unrestricted Frege systems).

\section{Algebraic formulation of $\fphp$ and a game with polynomials}

\label{3}

Let ${\bf F}_p[x_{i j}\ |\ i \in [n+1] \wedge j \in [n]]$
be the ring of polynomials over the finite field 
${\bf F}_p$ with $p$ elements with the indicated variables. 
Denote by 
$S_{n}$ the ring 
factored by the ideal generated by all polynomials $x_{i j}^2 - x_{i j}$.
Elements of $S_n$ are multi-linear polynomials.
Let $S_{n,e}$ be the ${\bf F}_p$-vector space of elements of
$S_n$ of degree at most $e$.
We shall denote monomials $x_a, \dots$ where $a, \dots$ are unordered tuples
of variable indices; the monomial is then the product of the
corresponding variables.

Beame et al.\cite{BIKPP} formulated (the negation of) PHP 
as the following {\bf $\fphp$-system}
of polynomial equations in $S_n$:

\begin{itemize}

\item $x_{i_1 j} \cdot x_{i_2 j} = 0$,  
for each $i_1 \neq i_2 \in [n+1]$ and $j \in [n]$.

\item $x_{i j_1} \cdot x_{i j_2} = 0$,  
for each $i \in [n+1]$ and $j_1 \neq j_2 \in [n]$.

\item $1 - \sum_{j \in [n]} x_{i j} = 0$, for each $i \in [n+1]$.

\end{itemize}
The left-hand sides of these equations will be denoted $Q_{i_1, i_2; j}$,
$Q_{i; j_1, j_2}$ and $Q_i$ respectively. 

The language of rings is a complete language for propositional logic
and it is easy to imagine a modification of the G-game to such a language
if the answers of Liar have to respect both the sum and the product.
The game we are going to define allows only simple questions and 
requires that sums of two polynomials and products of two
monomials are respected. 

We shall define the following game $H(e,n,r)$ played by two players Alice and Bob.
Alice's role will be similar to that of Prover in the G-game and Bob's to that of
Liar. In every round Alice may put to Bob a question of just one type:

\begin{enumerate}

\item[(A)] She asks Bob to give 
to a polynomial $f$ from $\sne$ a value from $\fp$.

\end{enumerate}
Bob's answers must obey the following rules:

\begin{enumerate}

\item[(B0)] If asked about a polynomial whose value he gave in an earlier
round Bob must answer identically as before.

\item[(B1)] He must give to each element $c \in \fp$ the value $c$,
and to each variable either $0$ or $1$.

\item[(B2)] If he gave values to $f$, $g$ and $f+g$, the values given to
$f$ and $g$ must sum up to the value he gave to $f + g$.

\item[(B3)] If he gave values to monomials  
$x_a$, $x_b$ and $x_a \cdot x_b$, the product of the
values given to $x_a$ and $x_b$ must equal to the value given to 
$x_a \cdot x_b$.

\item[(B4)] He must give value $0$ to all polynomials 
$Q_{i_1, i_2; j}$,
$Q_{i; j_1, j_2}$ and $Q_i$. 

\end{enumerate}
The game runs for $r$ rounds and Bob wins if he can answer all questions
while obeying the rules. Otherwise Alice wins.

We consider the multiplicativity condition for
monomials rather than for polynomials
as that more clearly isolates the role of linearity.
As is shown in Section \ref{4} the two versions
of the multiplicativity condition are essentially equivalent.

In principle Bob's strategy can be adaptive (i.e. his moves depend on
the development of the game) or even may depend on Alice. 
Call a strategy of Bob {\bf simple}
if it is a function $B$ assigning to elements of
$\sne$ values in $\fp$ and Bob, when asked to evaluate $f$, answers $B(f)$.
We shall abuse the language occasionally and talk about a {\bf simple Bob} rather than a
simple strategy for Bob.

\section{Five useful protocols for Alice}

\label{4}

In this section we describe five simple
protocols in which Alice can force Bob to answer various
more complicated questions,
similar to that of (P2). 

\bigskip

\noindent
{\bf Protocol $M_0$:} Assume that Bob asserted that $\sum_{i\le u}f_i \neq 0$.
Alice wants to force Bob to assert that $f_j \neq 0$ for some $j \le u$
(or to lose).

\medskip

She splits the sum into halves and asks Bob to evaluate
$\sum_{i\le u/2}f_i$ and $\sum_{i > u/2}f_i$.
As he already gave to  $\sum_{i\le u}f_i$ a non-zero value, 
by (B0) and (B2) - unless he quits - 
Bob must give to at least one of the half-sums a non-zero value.
Continuing in a binary search fashion in $\log u$ 
rounds she forces
Bob to assert that $f_j \neq 0$ for some $j \le u$.

\bigskip

\noindent
{\bf Protocol $M_1$:} Assume that Bob gave to some polynomials
$f$, $g$ and $f \cdot g$ values $B(f)$, $B(g)$ and
$B(f \cdot g)$ respectively, and that 
$B(f) \cdot B(g) \neq B(f \cdot g)$.
Alice wants to force Bob into a contradiction with the rules.

\medskip

Alice writes polynomials $f$ and $g$ as $\fp$-linear combinations
of monomials: $f = \sum_{a \in A} c_a x_a$ and
$g = \sum_{b \in B} d_b x_b$ with $c_a, d_b \in \fp$
and $x_a, x_b$ monomials. She
splits $A$ into two halves $A = A_0 \dot\cup A_1$, and asks Bob for the values
of 
$$
(\sum_{a \in A_0} c_a x_a)\ ,\ 
(\sum_{a \in A_0} c_a x_a)\cdot g\ ,\ 
(\sum_{a \in A_1} c_a x_a)\ ,\ \mbox{ and}\  
(\sum_{a \in A_1} c_a x_a)\cdot g
\ .
$$
Unless Bob violates the linearity rule (B2) his answers must satisfy
$$
B(\sum_{a \in A_i} c_a x_a)\cdot B(g)\ \neq\  
B((\sum_{a \in A_i} c_a x_a)\cdot g)
$$
for either $i = 0$ or $i = 1$.
Continuing in the binary search fashion Alice forces Bob to
assert
$$
B(c_a x_a)\cdot B(g)\ \neq\  
B( c_a x_a\cdot g)
$$
for some monomial $x_a$. Using (B1) and (B2) she forces
$$
B(c_a x_a) = c_a B(x_a)\ \mbox{ and }\  
B(c_a x_a g) = c_a B(x_a g) 
$$
and hence
$$
B(x_a) \cdot B(g) \ \neq\ B(x_a \cdot g)\ . 
$$
The number of variables is $n^{O(1)}$ and so the number
of monomials of degree at most $e$ is $n^{O(e)}$.
Hence all this process requires as most $O(e \log n)$
rounds of Alice's questions.

Now she analogously forces Bob to assert
$$
B(x_a) \cdot B(x_b)\ \neq \ B(x_a \cdot x_b)
$$
for some monomial $x_b$ occurring in $g$, violating thus (B3).

\bigskip

\noindent
{\bf Protocol $M_2$:} Assume that Bob asserted that $\Pi_{i\le k}f_i \neq 0$
and let $j \le k$ be arbitrary.
Alice wants to force Bob to assert that $f_j \neq 0$ (or to lose).

\medskip

She asks Bob to state the value of $f_j$ and if 
Bob says $f_j \neq 0$ she stops. Otherwise 
the triple $f_j, g$ and $f_j g$ for 
$g := \Pi_{i \le k, i\neq j}f_i$ satisfies the hypothesis
of protocol $M_1$  and Alice can win in $O(e \log n)$ rounds.

\bigskip

\noindent
{\bf Protocol $M_3$:} Assume that Bob asserted that $\Pi_{i\le k}f_i = 0$.
Alice wants to force Bob to assert that $f_j = 0$ for some $j \le k$ (or to 
lose).

\medskip

We shall describe the protocol by induction on $k$. Alice asks first 
for the value of $f_k$. If Bob states that $f_k = 0$ she stops. If he states that
$f_k \neq 0$ she 
asks him for the value of $\Pi_{i < k}f_i$.
If Bob says that $\Pi_{i < k}f_i = 0$, Alice 
has - by the induction hypothesis - a way how to solve the task. 

If he says that $\Pi_{i < k}f_i \neq 0$ 
Alice forces him into contradiction 
using protocol $M_1$. We may assume that all polynomials $f_i$
are non-constant and thus the induction process takes at most $k \le e$
steps.

Note that again Alice needed at most
$2e + O(e \log n) = O(e t)$  rounds in total.

\bigskip

\noindent
{\bf Protocol $M_4$:} Let $g = f^{p-1}$ and assume that
Bob gave to $g$ a value different from both $0, 1$.
Alice wants to force Bob into a contradiction.

\medskip

She asks Bob for the value of $f$ and assume Bob states $f = c \in \fp$.
If $c=0$ Alice uses protocol $M_2$ to force a contradiction.
If $c \neq 0$ Alice asks Bob for values of $f^2, f^3, \dots,
f^{p-1}$ and unless Bob returns values 
$c^2, c^3, \dots, c^{p-1}$ she forces him into a contradiction
by protocol $M_1$. But Bob cannot keep up these answers because
if he gave to $g$ now the value
$c^{p-1} = 1$ he would violate rule (B0).

\section{From Prover to Alice and from Bob to Liar}

\label{5}

In this section we employ the Razborov - Smolensky method
to show that the existence of many simple winning strategies 
for Bob yields a winning strategy for 
Liar\footnote{We could have bypassed the G-game and the explicit
use of the Razborov - Smolensky method by employing the characterization
of the size of $AC^0[p]$ Frege proofs in terms of degree of proofs in the
so called Extended Nullstellensatz of \cite{BIKPRS}. We prefer here a self-contained 
presentation.}. 
The reason to single out simple strategies is that we shall 
apply the Razborov - Smolensky approximation method 
in order to move from a G-game to an H-game,
by approximating formulas by low degree polynomials with respect to 
(a set of) Bob's strategies. 
The approximation process (and hence a strategy for Alice 
to be constructed) depends on the set of Bob's strategies
we start with 
and to avoid circularity we restrict to sets containing only (but not
necessarily all) simple strategies.

\begin{lemma} \label{5.1}
Let $d \geq 2$, $n \geq 1$ and $t \geq \log n$ be arbitrary 
and take
parameters $e, r$
$$
e\ :=\ ((t^2 + 2t) p)^{d}  \ \mbox{ and }\ 
r\ :=\ O(e t^{4})\ .
$$
Let $P$ be any strategy for Prover in game $G(d, n, t)$.
Let $\on$ be a non-empty set of simple strategies for Bob in game $H(e, n, r)$.

Then $P$ can be translated into a strategy $A$ for Alice in $H(e, n, r)$ such that 
the following holds:
\begin{itemize}

\item
If
\begin{equation}\label{e}
\prob_{B \in \on}[B\ \mbox{ wins over $A$ in } H(e, n, r)]\ >\ 
1\ -\ 2^{-(t +1)}
\end{equation}
then there exists a Liar's strategy $L$ winning over $P$ in $G(d, n, t)$.
\end{itemize}
\end{lemma}

\prf

Let $P$ and $\on$ be given. Let $F$ be the smallest set of formulas
closed under subformulas and containing all possible $P$'s questions 
according to rule (P1) in all
plays of the game $G(d, n, t)$ against all possible Liars. 
The number of such (P1) questions is at most $2^{t^2}$
and each has size at most $2^t$ and so also at most $2^t$
subformulas.
Thus 
the depth of all formulas in $F$ is at most $d$ and their total number is
bounded by $2^{t^2 + t}$.

We shall use the Razborov - Smolensky method to assign to all formulas
$\varphi \in F$ a polynomial $\hat \varphi \in \sne$. However, we shall
approximate with respect to Bob's 
strategies from $\on$ rather than with respect to all
assignments to variables as it is usual.

Fix parameter $\ell := t^2 + 2t$. Put
$\hat x_{e} := x_{e}$, $\hat{(\neg \varphi)} := 1 - \hat \varphi$ and for
$\varphi = MOD_{p,i}(\varphi_1, \dots, \varphi_k)$ define
$$
\hat{\varphi} \ :=\  1\ -\ 
((\sum_{j \le i} \hat \varphi_j) -i)^{p-1}\ .
$$
For the remaining case $\varphi = \bigvee_{i \in [u]}\varphi_i$
assume that all polynomials $\hat \varphi_i$ were already defined.
Pick $\ell$ subsets
$J_1, \dots, J_{\ell} \subseteq [u]$, 
independently and uniformly at random (we shall fix them in a moment), 
and define polynomial
$$
p_{\varphi}(y_1, \dots, y_u)\ :=\ 
1\ -\ \Pi_{j \le \ell} (1\ -\ (\sum_{i \in J_j} y_i)^{p-1})\ 
$$
and using $p_{\varphi}$ put
\begin{equation} \label{approx}
\hat \varphi\ :=\ p_{\varphi}(\hat \varphi_1, \dots, \hat \varphi_u)\ .
\end{equation}
The following claim is easily verified by induction on the depth of $\varphi$, using
the protocols from Section \ref{4}.

\medskip

\noindent
{\bf Claim 1:} {\em
Let $\varphi \in F$ and assume that Bob asserted that $\hat \varphi = c \in \fp$ for some
$c \neq 0, 1$. Then Alice can force Bob into a contradiction in
$O(e\log n)$ rounds.
}

\medskip

Let $b_i \in \{0,1\}$ be the truth-value of the statement $B(\hat \varphi_i)\neq 0$.
For $B \in \on$ we have that
\begin{equation}\label{rs}
\bigvee_{i \in [u]} b_i\ =\ 
p_{\varphi}(b_1, \dots, b_u)
\end{equation}
with the probability at least $1 - 2^{- \ell}$ (taken over the choices of sets $J$). 
Hence we can select specific 
sets $J_1, \dots, J_{\ell}$ such that (\ref{rs}) holds for all but 
$2^{-\ell}\cdot |\on|$ simple Bobs from $\on$.
The polynomial $\hat \varphi$ in (\ref{approx}) is assumed to
have this property.

\medskip

Define in this way the polynomial $\hat \varphi$ for
all (at most $2^{t^2 + t}$)
formulas $\varphi \in F$ by induction on the depth $1, 2, \dots, d$.
Each is of degree at most $(\ell (p-1))^d \le ((t^2 + 2t)p)^d = e$ and 
it holds that:

\medskip
\noindent 
{\bf Claim 2:} {\em There is a subset $Err \subseteq \on$ such that
$|Err| \le 2^{- t}|\on|$ and such that (\ref{rs}) holds for all
disjunctions $\varphi \in F$ and all $B \in \on \setminus Err$.}
\medskip

Now we define, using the given
strategy $P$ for Prover, a specific strategy $A$ for
Alice in $H(e, n, r)$. We transcript $P$ into $A$ a question 
by question; each question of $P$ may be replaced by a series of questions of
Alice.

If P asks according to (P1) what is the value of $\varphi$,
Alice simply asks for the value of $\hat \varphi$.
Let $\varphi = \bigvee_{i \in [u]}\varphi_i$ and assume that
P asks according to (P2); there are two cases to consider: 

\begin{enumerate}
\item[(a)] $\varphi$ got value {\sf false} and P asks for the value of 
one disjunct $\varphi_j$,

\item[(b)] $\varphi$ got value {\sf true} and P asks for a witness $\varphi_j$.
\end{enumerate}

Assume for the case (a) that Bob asserted in an earlier round that
$\hat\varphi = 0$. Alice asks Bob for the value of $\hat \varphi_j$. If he gives
$B(\hat \varphi_j) = 0$ the simulation of P  moves to the next round. If he 
replies that $B(\hat \varphi_j) = 1$,
Alice uses first protocol
$M_2$ repeatedly to force Bob to assert
$$
1\ - \ (\sum_{i \in J_v} \hat\varphi_i)^{p-1}\ \neq \ 0
$$
for all $v \le \ell$.
Then for each $v$ she uses protocol $M_4$ to force Bob to say that 
$$
(\sum_{i \in J_v} \hat\varphi_i)^{p-1}\ = \ 0
$$
and further protocol $M_3$ to assert that
\begin{equation}\label{a}
\sum_{i \in J_v} \hat\varphi_i\ = \ 0\ .
\end{equation}
This needs $O(\ell e \log n) = O(t^2 e \log n) = O(e t^3)$
rounds. 

As $B(\hat \varphi_j) = 1$, if
Bob uses a strategy $B \in \on \setminus Err$,
the definition of $Err$ guarantees that one of the equations 
in (\ref{a}) is false when $\hat \varphi_i$'s are 
evaluated by $B$:
$$
B(\sum_{i \in J_v} \hat\varphi_i)\ \neq\ 
\sum_{i \in J_v} B(\hat\varphi_i)\ .
$$
This itself is not a violation of rule (B2) but Alice can use this situation
and to force Bob to lose. We shall describe her strategy
as probabilistic; a deterministic one is obtained by an averaging argument.

Alice splits $J_v = K_0 \dot \cup K_1$ into halves
and asks Bob for the values of $\sum_{i \in K_0} \hat\varphi_i$
and $\sum_{i \in K_1} \hat\varphi_i$. Unless he violates (B2)
his answers must sum up to $B(\sum_{i \in J_v} \hat\varphi_i)$.
Hence for $k=0$ or $k=1$
$$
B(\sum_{i \in K_k} \hat\varphi_i)\ \neq\ 
\sum_{i \in K_k} B(\hat\varphi_i)\ .
$$ 
Alice guesses for which $k$ this happens and then proceeds analogously
with $\sum_{i \in K_k} \hat\varphi_i$, splitting it into halves,
asking Bob for the values, etc. If she always guesses right
then in $t$ steps (as the size of the sums is bounded by $2^t$)
she will reduce the sums to one term and will win.
Alice has the probability at least $2^{-t}$ to make the right 
choices. She does not know a priori which of the $\ell$
sums $\sum_{i \in J_v} \hat\varphi_i$ to use so she must try all.
This takes $O(\ell t) = O(t^3)$ rounds.

There are at most $t$ simulations of a (P2a) question in the G-game
but Alice needs to employ the random strategy above only 
once when the case $B(\hat \varphi_j) = 1$ occurs, and then her probability
of success is at least $2^{-t}$.
By averaging there are fixed choices that Alice can make, yielding this
success probability outside of $Err$. In particular, for a random 
$B \in \on \setminus Err$, if Alice uses these choices then
either B must give to $\hat \varphi_j$ value $0$ or Alice wins
with the probability at least $2^{-t}$.
We shall describe this situation below by the phrase
that the {\em (P2a) simulation succeeded}.

\medskip

Assume for the case (b) that Bob answered earlier 
that $\hat \varphi = 1$
and hence also that
$$
\Pi_{j \le \ell}(1\ -\ (\sum_{i \in J_j}\hat\varphi_i)^{p-1})\ =\ 0.
$$
Alice uses protocols $M_3$ and $M_4$
to force Bob to state that
$\sum_{i \in J_v}\hat \varphi_i = 1$ for some $v \le \ell$.
This uses $O(e \log n) = O(e t)$ rounds.
Then she uses protocol $M_0$
to force Bob to say that
$\hat \varphi_j \neq 0$ for some $j \in J_v$
and by Claim 1 the value has to be $1$ ($O(e \log n) =
O(e t)$ rounds are used in Claim 1).
The number of formulas $\varphi_i$ is bounded by the size of
$\varphi$, i.e. by $2^t$, and so this uses at most $t$ rounds in protocol $M_0$,
i.e. still $O(e t)$ in total.

This describes the strategy $A$. 

\bigskip

By Claim 2 with the probability at least $1 - 2^{-t}$
a random $B \in \on$ is outside $Err$, and for these Alice's
simulations of (P2a) questions succeed with the probability at least
$2^{-t}$. Thus the
inequality (\ref{e}) from the
hypothesis of the lemma implies that there is at least one
$B \in \on \setminus Err$ winning over the particular Alice's strategy
A and for which A's  simulations of (P2a) questions succeed.

Use B to define a strategy L for Liar in the original game $G(d, n, t)$
simply by giving to $\varphi$ the truth value 
$B(\hat\varphi)$ when asked a (P1) type question, and
giving a witness $\varphi_j$ constructed in the case (b) above
when asked a (P2b) type question.

From the construction of A (and rules for Bob) it follows that L
satisfies the rules for Liar. In particular, by (B4) all polynomials 
from the $\fphp$-system
 get $0$ by B and so all axioms of $\fphp$ 
get by L value
{\sf true}.

Note that one question of P is transcribed into at most
$O(e t^3)$ Alice's questions.
Hence in every play of the H-game transcribing a play of the $G$-game
there are in total  at most 
$r  = O(e t^4)$ rounds. 

\qed

\section{A general reduction to a search problem}
\label{6}

The reduction of the lengths-of-proofs problem to a question
about the $H$-games in Sections \ref{2} - \ref{5} is not specific to
$\fphp$ and works in a fairly general situation that we shall describe now.
Then we reduce the proof complexity problem further to a question
about the computational complexity of a certain task involving computations
with search trees.

The only specific thing in the $\fphp$ case is 
the transcription of the axioms of $\fphp$ into the $\fphp$
polynomial system in Section \ref{3}. This is not a mere mechanical translation
from DeMorgan language into the language of rings (in that the 
axioms $\bigvee_{j \in [n]} p_{i j}$
 would translate into polynomials of degree about $n$
and not into degree $1$ polynomials $Q_i$). 
In order to avoid inevitable technicalities when trying to define
suitable translations from a general set of axioms to 
a polynomial system we simply take
as our starting point an unsolvable system of polynomial equations
of a constant degree. The 
truth value of an equation $f(x_1, \dots, x_m)=0$ for Boolean variables $x_i$,
$f$ a degree $O(1)$ polynomial over $\fp$, can be defined by a
depth $2$, size $m^{O(1)}$ 
$AC^0[p]$ formula. Namely, writing $f$ as an $\fp$-linear combination
$\sum_{a \in A} c_a x_a$ of monomials $x_a$, with $c_a \in \{1, \dots, p-1\}$,
consider the formula
\begin{equation}
\varphi \ :=\ 
MOD_{p,0}(\psi_1, \dots, \psi_k)
\end{equation}
where $k = \sum_{a \in A} c_a$ and $\psi_i$'s are
conjunctions of variables corresponding
to monomials from $f$, each monomial $x_a$ being represented 
$c_a$-times. Clearly\footnote{Instead of assuming degree $O(1)$ it would suffice
to assume that $f$ is an $\fp$-linear combination of polynomially
many monomials.} $\varphi$ represents the truth value of $f=0$ on
Boolean variables.
The polynomial system can thus be also thought of
as an unsatisfiable set of $AC^0[p]$ formulas and we can speak about
its $LK_d(MOD_p)$ refutations.

\bigskip

We shall now consider the following general set-up. For $n = 1, 2, \dots$
let ${\cal F}_n$ be a sequence of sets of polynomials over $\fp$ 
in variables $Var({\cal F}_n)$. We shall
assume that: 
\begin{enumerate}
\item 
polynomials in sets ${\cal F}_n$ have $O(1)$ degree,

\item the size of both ${\cal F}_n$ and 
$Var({\cal F}_n)$ is $n^{O(1)}$,

\item the polynomial system
$$
f\ =\ 0\ ,\ \mbox{ for } \ f \in {\cal F}_n
$$
contains equations $x^2 - x = 0$ for all 
$x \in Var({\cal F}_n)$
and is unsolvable in $\fp$.
\end{enumerate}
Let $\sfn$ be the $\fp$ - vector space of multi-linear polynomials
in variables of ${\cal F}_n$ and of degree at most $e$.

We want to replace games and strategies considered in previous sections by
a more direct computational model, namely that of search trees. 
Define an {\bf $\sfn$ - search tree} $T$ to be a $p$-ary tree whose
inner nodes (non-leaves) are labelled by polynomials from
$\sfn$, the $p$ edges leaving a node labelled by $g$
are labelled by $g=0, g=1, \dots, g = p-1$, and leaves
are labelled by elements of a set $X$. 

Any function $B : \sfn \rightarrow \fp$ determines a path 
$P_T(B)$ in $T$ 
consisting of edges labelled by $g = B(g)$
and thus it also determines an element of $X$: the label of the
unique leaf on $P_T(B)$. Hence $T$ defines a function assigning
to any map $B : \sfn \rightarrow \fp$ an element of $X$ to be denoted $T(B)$.

Let $\eer$ be the set of pairs and triples of the form
$(B1, c)$ for $c \in \fp$ or $(B1, x)$ for $x \in Var({\cal F}_n)$,
$(B2, f, g)$, $(B3, x_a, x_b)$ or $(B4, f)$ for $f \in {\cal F}_n$,
with $f,g, x_a, x_b$ of degree at most $e$. These are intended to
indicate what instance of which
rule did Bob violate. We say that
$(B1, c)$ is an {\bf error for} $B$ iff $B(c) \neq c$, $(B1, x)$
is an error for $B$ iff $B(x) \neq 0, 1$, and similarly for
the other pairs and triples\footnote{We ignore errors for (B0) as that rule cannot 
be violated by a simple Bob
and hence search trees do not need to ask anything twice on any path.}.

In the following statement we talk about refutations of equations $f=0$, 
$f \in {\cal F}_n$. As pointed out earlier, we can view them also as depth $2$,
polynomial size formulas with $MOD_{p,0}$ connectives and hence it makes
a prefect sense to talk about their $LK_d(MOD_p)$-refutations.

The reductions of Sections \ref{2} - \ref{5} 
used the example of $\fphp$ (see the beginning of Section \ref{2}) but nothing
specific to it was used. Hence we can employ 
the reductions to derive the following general statement. 
In it we replace degree $e$ by (bigger) $r$
in order to avoid the need to define here
the relation between them implicit 
in Lemma \ref{5.1}.

\begin{theorem} \label{6.2}
Let $r = r(n) \geq (\log n)^{\omega(1)}$ be a function and 
let ${\cal F}_n$ be sets of polynomials obeying the restrictions
1., 2. and 3. listed above. 

Then for every $d \geq 2$ there are $\epsilon_d > 0$
and $n_d \geq 1$ such that for an arbitrary non-empty
set $\om$ of maps from $\sfnr$ to $\fp$
the following 
implication (I) holds for all $n \geq n_d$ and all 
$0 < \epsilon \le \epsilon_d$:

\begin{itemize}
\item[(I)]
If for every $\sfnr$ - search tree $T$
of depth $r$ and with
leaves labelled by elements of $\eerr$
it holds that 
\begin{equation} \label{key}
\prob_{B \in \om}[\mbox{ $T(B)$ is not an error for $B$ }]
\ > \ 1 -
2^{- r^{\epsilon}}
\end{equation}
then 
$LK_d(MOD_p)$ does not refute the set of
formulas $f=0$, $f \in {\cal F}_n$, by a proof of size
less than $2^{\Omega(r^{\epsilon})}$.
\end{itemize}

\end{theorem}

\prf

Assume that $LK_d(MOD_p)$ does refute the set of
formulas $f=0$, $f \in {\cal F}_n$, by a proof of size
$s = s(n)$. By Lemma \ref{2.1} Prover has a winning strategy P
for game $G(d+c, n, t)$, where $t = t(n) = O(\log s)$ and $c$ is
an absolute constant.

Put $\epsilon_d := \frac{1}{2(d+c)+5}$ and let $0 < \epsilon \le
\epsilon_d$. If it were that $t+1 \le r^\epsilon$ then  
the parameters $e', r'$ of the game $H(e',n,r')$ constructed
in Lemma \ref{5.1} 
satisfy 
$e' \le r' < r$ and, in particular, the game is an $H(r,n,r)$ game.

The strategy A defined in Lemma \ref{5.1} for the game
defines an  $\sfnr$ - search tree $T$
of depth $r$ and with
leaves labelled by elements of $\eerr$ in a natural way:
a path in $T$ corresponds to possible answers of a simple Bob
and the path stops as soon as a violation of one of the rules (B1)-(B4) 
occurs (rule (B0) cannot be broken by a simple Bob).
The label of the resulting leaf is the instance of the rule
that was broken
(if a violation did not occur we use any element of $\eerr$).

Assume that $\om$ is a set of simple Bobs
for which the inequality (\ref{key}) holds. Then also 
the inequality (\ref{e}) from Lemma \ref{5.1} holds and thus
by that lemma there is a strategy L for Liar that wins over P in
the original $G$-game. That is a contradiction and thus
$s \geq 2^{\Omega(r^{\epsilon})}$.

\qed

\bigskip

To conclude the paper let us discuss informally the construction
underlying Lemma \ref{5.1} and Theorem \ref{6.2}. In particular, we
see these formal statements as templates for a possible variety
of analogous reductions, and it is not clear which one - if any - will
be eventually useful. 

The strategy $A$ is constructed in Lemma \ref{5.1} by a randomized
process from strategy $P$ and from set $\on$. Let us call the class
of all strategies
$A$ that can occur in this way the class of {\em $(P, \on)$-generated
strategies}. One such class contains only a few of 
all possible Alice's strategies.
Moreover, we can pick $\on$ depending on $P$. Hence one can weaken the
hypothesis in these statements and, for example, Theorem \ref{6.2} 
could be reformulated as follows:  

\begin{itemize}

\item
Let $r = r(n) \geq (\log n)^{\omega(1)}$ be a function and 
let ${\cal F}_n$ be sets of polynomials obeying the restrictions
1., 2. and 3. listed above. 

Then for every $d \geq 2$ there are $\epsilon_d > 0$
and $n_d \geq 1$ such that the following holds:

If for every Prover's strategy $P$ for game $G(d,n,r^{\Omega(1)})$
there exists a non-empty
set $\om(P)$ of maps from $\sfnr$ to $\fp$
then the following 
implication (I') holds for all $n \geq n_d$ and all 
$0 < \epsilon \le \epsilon_d$:

\begin{itemize}
\item[(I')]
If for every $\sfnr$ - search tree $T$
of depth $r$ and with
leaves labelled by elements of $\eerr$
originating from a $(P, \om(P))$-generated
$A$
it holds that 
\begin{equation}
\prob_{B \in \om(P)}[\mbox{ $T(B)$ is not an error for $B$ }]
\ > \ 1 -
2^{- r^{\epsilon}}
\end{equation}
then 
$LK_d(MOD_p)$ does not refute the set of
formulas $f=0$, $f \in {\cal F}_n$, by a proof of size
less than $2^{\Omega(r^{\epsilon})}$.
\end{itemize}

\end{itemize}

\medskip
\noindent
This formulation stains the combinatorially clean original
formulation by a reference to $P$ but (I') may be a weaker hypothesis
to arrange.

Another issue is the discouragingly high probability required in 
(\ref{e}) and (\ref{key}). This is due solely by Alice's simulation
of the (P2a) move of $P$. At that point she found $\ell \le
O(t^2) \le e \le r$ sets $K$, $|K| \le 2^t$, of degree $e$ polynomials
such that for one of them $B$ fails linearity:
\begin{equation}\label{concl}
B(\sum_{i \in K} g_i)\ \neq\ 
\sum_{i \in K} B(g_i)\ 
\end{equation}
and her strategy worked up to this point
for all $B \notin Err$ (as long as $P$
was a winning strategy for the Prover). Getting from this situation to
a violation of rule (B2) costs her the drop of the success probability by 
the multiplicative factor $2^{-t}$. Hence we could redefine the rules
for the H-game and, in particular, 
the error sets $\eer$ for the search problems
to be solved by the trees, and include that situation  
(i.e. $A$ producing $\ell$ sets $K$ such that one of them satisfies
(\ref{concl})) among the stopping Bob's errors. Let us call 
$*\eer$ the set of errors with this new type of an error added.
Then we could reformulate Theorem \ref{6.2} 
differently as follows:  

\begin{itemize}

\item
Let $r = r(n) \geq (\log n)^{\omega(1)}$ be a function and 
let ${\cal F}_n$ be sets of polynomials obeying the restrictions
1., 2. and 3. listed above. 

Then for every $d \geq 2$ there are $\epsilon_d > 0$
and $n_d \geq 1$ such that for an arbitrary non-empty
set $\om$ of maps from $\sfnr$ to $\fp$
the following 
implication (I'') holds for all $n \geq n_d$ and all 
$0 < \epsilon \le \epsilon_d$:

\begin{itemize}
\item[(I'')]
If for every $\sfnr$ - search tree $T$
of depth $r$ and with
leaves labelled by elements of $*\eerr$
it holds that 
\begin{equation}
\prob_{B \in \om}[\mbox{ $T(B)$ is not an error for $B$ }]
\ > \ 
2^{- r^{\epsilon}}
\end{equation}
then 
$LK_d(MOD_p)$ does not refute the set of
formulas $f=0$, $f \in {\cal F}_n$, by a proof of size
less than $2^{\Omega(r^{\epsilon})}$.
\end{itemize}
\end{itemize}

\medskip
\noindent
Let us stress that the culprit property is the linearity
by observing that simple Bobs can be without a loss of generality assumed
to satisfy all rules except possibly (B2). First, having $B$ we can
define $B'$ by correcting all values of $B$ that violate rules (B1) or
(B4). If $B'$ is asked by Alice for one of these new values, the original $B$
would lose. Hence $B'$ is as good as $B$ against any $A$.

Then define $B''$ by giving to every monomial $x_a = \Pi_i x_i$
the value $\Pi_i B'(x_i)$. Enhance any $A$ to vigilant $A^*$
that whenever she asks for the value of a monomial, she asks also for the
values of all its variables (this enlarges the number of round $e$-times at 
most). Clearly, $B''$ fares as well as $B'$ against a vigilant $A^*$.

Finally, let us remark that
it would be interesting and possibly quite useful to
modify the construction so that 
adaptive Bobs are allowed.

\bigskip
\noindent
{\bf Acknowledgements.}

I thank L.~Kolodziejczyk (Warsaw),
S.~M\"{u}ller (Tokyo), P.~Pudl\' ak (Prague) and N.~Thapen (Prague)
for critical comments and discussions.
I am also indebted to the anonymous referee for helpful suggestions.

\bigskip
\noindent
{\bf Mailing address:}

Department of Algebra

Faculty of Mathematics and Physics

Charles University

Sokolovsk\' a 83, Prague 8, CZ - 186 75

The Czech Republic

{\tt krajicek@karlin.mff.cuni.cz}

\end{document}